# HYPERSPECTRAL CLASSIFICATION WITH ADAPTIVELY WEIGHTED L₁-NORM REGULARIZATION AND SPATIAL POSTPROCESSING

*Victor Stefan Aldea*
victor.aldea@mail.mcgill.ca

**ABSTRACT**

**Sparse regression methods have been proven effective in a wide range of signal processing problems such as image compression, speech coding, channel equalization, linear regression and classification. In this paper a new convex method of hyperspectral image classification is developed based on the sparse unmixing algorithm SUnSAL for which a pixel adaptive $L_1$-norm regularization term is introduced. To further enhance class separability, the algorithm is kernelized using a RBF kernel and the final results are improved by a combination of spatial pre and post-processing operations. It is shown that the proposed method is competitive with state of the art algorithms such as SVM-CK, KSOMP-CK and KSSP-CK.**

## I. INTRODUCTION

Over the last fifteen years hyperspectral images (HSI), which record the electromagnetic spectrum in a few tens to a few thousands of spectral bands, have been used in a variety of tasks such as hyperspectral image classification [1], [2], [3], target detection [5], [6], analytical chemistry, astronomy, pharmaceutical process modelling and biomedical applications [7].

A few overview papers written on this subject are [7], [8], [9], [10], together with the March 2013 issue of the Proceedings of the IEEE and the January 2014 issue of the IEEE Signal Processing Magazine. For image classification purposes, the available methods can generally be subdivided into statistical, neural networks and sparse regression-based algorithms.

Traditionally, many supervised statistical classifiers use the training data set to build models of the underlying density in the feature space for each of the various classes in the training set (mostly by using a mixture of Gaussians) [11], [12], [13]. Besides the fact that the Gaussian assumption is often incorrect, this density estimation in high-dimensional spaces suffers from the Hughes effect: for a fixed amount of training data, the classification accuracy as a function of the dimension of the data (the number of hyperspectral bands) reaches a maximum and then declines with increasing dimension, because there is a limited amount of training data to estimate the larger and larger number of parameters needed to describe the model of the densities of each class. To deal with this, usually a feature selection/reduction step is first performed on the high-dimensional data to reduce its dimensionality. Thus statistical classifiers either require considerable pre-processing or may entirely fail to work correctly on high-dimensional data.

Another well-established method of supervised classification is the SVM (support vector machine) [11 ch.3], [14], which separates classes by a hyperplane whose parameters are established using the training data (dictionary). While it has proven itself immune to the Hughes effect, it needs to be retrained for every single change to its dictionary and, in order to achieve good classification results it needs to use spatial-spectral composite kernels [15]. However, the spatial information required here may not always be available as part of the provided dictionary. Other variations of the SVM have been proposed to further improve classification performance when the number of training samples is low such as the transductive SVM [16], [17], [18] which increases its original training set by iteratively including initially unlabeled samples in a semi-supervised manner, or LFDA-SVM [19].

Due to the success of sparse coding for image compression and face recognition applications [20], many sparse regression (SR) algorithms have been proposed for the domain of image and signal processing [21], [22], [23], including the hyperspectral pixel unmixing problem. Unlike the SVM which can be considered to be a discriminative method, SR algorithms can be seen as generative models, where the subspaces representing different classes compete with each other during the pixel unmixing process, leading to a vector of unmixing coefficients which has only a few non-zero, representative coefficients.

Recently, a fast and efficient method for spectral unmixing, the SUnSAL algorithm, has been proposed in [24], [25]. It has already been used with very good results for HSI classification in several papers [3], [4] and generalized to a spatial-spectral collaborative method, CL-SUnSAL, in [26], [27].
Given that the CL-SUnSAL method is slower converging than SUnSAL and not directly kernelizable, an adaptively weighted $L_1$-norm SUnSAL algorithm is proposed together with spatial post-processing in order to achieve competitive classification results.

The rest of this paper is organized as follows: Section II provides the necessary background concerning sparse unmixing via the SUnSAL algorithm, Section III introduces the development of the adaptively weighted $L_1$-norm SUnSal algoritm with spatial postprocessing for classification and Section IV summarizes our experimental results. Concluding remarks are presented in Section V.

## II. BACKGROUND

Given the dictionary matrix $A \in R^{k \times n}$ and the observed mixed hyperspectral pixel $y \in R^k$, let $x \in R^n$ be the vector of unknown unmixing coefficients of the columns of A. For classification purposes, there are several signatures (i.e. columns) of $A$ for each class, and therefore usually $n \gg k$. Since it has been observed that a hyperspectral pixel is made up usually of a reduced number of spectral signatures (endmembers) compared to the total number of endmembers present in a given image, [e.g. 7], we know that the solution vector $x$ is sparse, i.e. it contains only a few non-zero entries, and should be obtained as the result of a constrained sparse regression (CSR) problem such as the $L_1$-norm regularized optimization:

$$P_{CSR}: \quad \min_{x} \frac{1}{2}\|Ax - y\|_2^2 + \lambda \|x\|_1 \quad (1)$$

possibly subject to $x \geq 0$, where $\lambda$ is the parameter controlling the relative weight between the $L_2$ and $L_1$ terms.

To solve problems of the type shown in eq.(1) it has been observed that it is often easier and more efficient to solve problems of the type shown in eq.(2) via an Augmented Lagrangian method in which variable splitting has been introduced. To this effect, the ADMM (alternating direction method of multipliers) has been developed in [28] and further specialized in [24], [25]. It can briefly be described as a general optimization problem of the type:

$$\min_{x \in R^n} f_1(x) + f_2(Gx) = \min_{x \in R^n, u \in R^p} f_1(x) + f_2(u) \quad (2)$$

where $f_1: R^n \to \bar{R}$, $f_2: R^p \to \bar{R}$ and $G \in R^{p \times n}$, $\bar{R} = R \cup \infty$. A summary of the ADMM steps is shown in **Algorithm** 1.

**Algorithm 1.** ADMM algorithm [25, 28]
1. Set k=0, choose $\mu > 0$, $u_0$ and $d_0$.
2. **Repeat**

$$3. x_{k+1} \in \arg\min_{x} f_1(x) + \frac{\mu}{2}\|Gx - u_k - d_k\|_2^2.$$

$$4. u_{k+1} \in \arg\min_{u} f_2(u) + \frac{\mu}{2}\|Gx_{k+1} - u - d_k\|_2^2.$$

5. $d_{k+1} \leftarrow d_k - (Gx_{k+1} - u_{k+1})$.
6. $k \leftarrow k+1$.
7. **Until** stopping criterion is satisfied.

## III. PROPOSED WEIGHTED $L_1$-NORM SUnSAL

While $L_1$-norm regularized cost functions such as (1) are convex and provide some measure of sparsity, they do not always provide sufficient sparsity and result in an incorrect recovery of the support of the signal $y$ in terms of the columns of the dictionary matrix $A$, which is a structural error [29]. On the other hand, it is known that $L_p$-norm regularized cost functions, $0 \leq p < 1$, are much more effective at recovering the correct sparse support of the signal $y$ as $p$ approaches zero, but their solutions are increasingly prone to be caught in some local minimum as $p$ decreases (convergence error), since they are not convex functions [29]. Many efforts have been deployed to obtain convex approximations to $L_p$-norm regularized cost functions with increasing success [30], [31], [32], [33].

To address this shortcoming of $L_1$-norm regularized cost functions while still working in the context of a convex optimization problem, we introduce a diagonal weighting matrix $\Gamma$ meant to further enhance the sparsity of the solution:

$$P_{CSR}: \quad \min_{x} \frac{1}{2}\|Ax - y\|_2^2 + \lambda \|\Gamma x\|_1 \quad (3)$$

A diagonal entry in $\Gamma$ is large if the corresponding column of $A$ is far from the current hyperspectral pixel $y$ and small if the two vectors are close together, thus representing how much importance the various columns of $A$ and their corresponding regression coefficients in $x$ should be given in the optimization problem[1]. Many measures of closeness can potentially be used for the diagonal entries of $\Gamma$, such as the Euclidean distance between $y$ and each column of $A$, the angle, or the Euclidean distance after using LFDA [35], [36].

The vector of weights on the diagonal of $\Gamma$ is then linearly re-scaled to an appropriate range= [*minWeight, maxWeight*]; and afterwards non-linearly transformed by the application of the *tanh()* function. This procedure is repeated two or three times to ensure that the number of entries on the diagonal of $\Gamma$ still having small values has become now very small compared to the initial distribution, as in Fig.1. The number of iterations of the above procedure and the range are usually set by validation, but several experiments have shown that the range= [1.42, 3.50]; and 2-3 iterations usually work well in most cases.

---
[1] Note that introducing $\Gamma$ in eq.(3) instead of just using the $L_1$-norm of $x$ could be seen as using an adaptive dictionary $A$ for the unmixing of each hyperspectral pixel $y$.

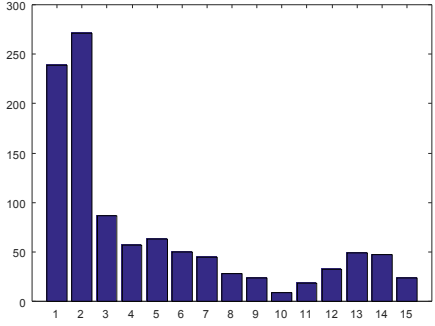

a) Original Weight distribution in $\Gamma$ after 0 iterations.

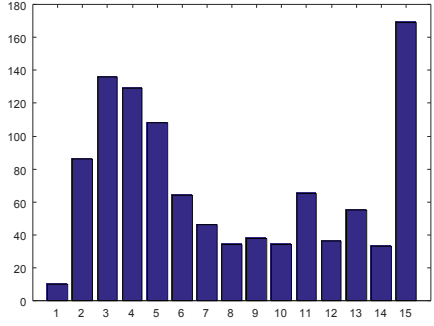

b) Modified Weight distribution in $\Gamma$ after 1 iteration.

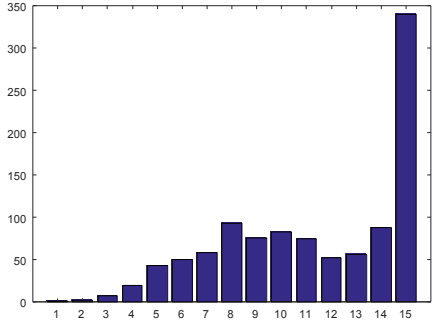

c) Modified Weight distribution in $\Gamma$ after 2 iterations.

**Fig.1** Succesive modification of the weight distribution on the diagonal of the adaptive $\Gamma$ matrix.

Using the following definitions for $f_1(x)$ and $f_2(u)$ in eqs.(2) and (3)

$$f_1(x) = \frac{1}{2}\|Ax - y\|_2^2$$
$$f_2(u) = \lambda \|u\|_1 \quad (4)$$
$$G = \Gamma, \; u = \Gamma x.$$

Step 3. of the ADMM procedure requires solving a quadratic problem

$$x_{k+1} \in \arg\min_{x} f_1(x) + \frac{\mu}{2}\|\Gamma x - u_k - d_k\|_2^2 =$$
$$\arg\min_{x} \frac{1}{2}\|Ax - y\|_2^2 + \frac{\mu}{2}\|\Gamma x - u_k - d_k\|_2^2. \quad (5)$$

with Hessian and linear term given by

$$H = A^T A + \mu \Gamma^T \Gamma$$
$$f = A^T y + \mu \Gamma^T (u_k + d_k) \quad (6)$$

and solution $x_{k+1} = H^{-1} f$.

Step 4. of the ADMM procedure becomes

$$u_{k+1} \in \arg\min_{u} f_2(u) + \frac{\mu}{2}\|\Gamma x_{k+1} - u - d_k\|_2^2 =$$
$$\arg\min_{u} \frac{\mu}{2}\|\Gamma x_{k+1} - u - d_k\|_2^2 + \lambda \|u\|_1 = \quad (7)$$
$$\arg\min_{u} \frac{1}{2}\|u - v_k\|_2^2 + \frac{\lambda}{\mu}\|u\|_1.$$

where we have used $v_k = \Gamma x_{k+1} - d_k$.

The solution to (7) is the soft-threshold function [25, 37]

$$u_{k+1} = soft\left(v_k, \frac{\lambda}{\mu}\right) \quad (8)$$

To impose the positivity constraint on the solution one only needs to project the result in (8) onto the positive orthant:

$$u_{k+1} = \max\left\{0, \; soft\left(v_k, \frac{\lambda}{\mu}\right)\right\} \quad (9)$$

A summary of the SUnSAL algorithm steps with the adaptively weighted $L_1$-norm is shown in **Algorithm** 2.

**Algorithm 2**. Weighted $L_1$-norm SUnSAL algorithm with pixel adaptive regularization (Proposed Algorithm).
1. Set k=0, choose $\mu > 0$, $u_0$ and $d_0$.
2. **Repeat**
   3. $H = A^T A + \mu \Gamma^T \Gamma, \; f_k = A^T y + \mu \Gamma^T (u_k + d_k)$
   4. $x_{k+1} = H^{-1} f_k$.
   5. $v_{k+1} = \Gamma x_{k+1} - d_k$.
   6. $u_{k+1} = soft\left(v_{k+1}, \frac{\lambda}{\mu}\right)$.
   7. $d_{k+1} \leftarrow d_k - (\Gamma x_{k+1} - u_{k+1})$.
   8. $k \leftarrow k+1$.
9. **Until** stopping criterion is satisfied.

Unlike the original SUnSAL [25] algorithm, the adaptively weighted $L_1$-norm variant proposed in this section introduces an adaptive Hessian and linear term for each pixel which is unmixed. Since the two terms in (3) are closed, proper, convex functions[2] and $\Gamma$ is designed to always be full-column rank, the procedure in **Algorithm 2** always converges (see Thm.1 from [25], [28]). Because the singular value decomposition used to invert the Hessian in the original paper [25] is too computationally complex to be used for every iteration, it has been replaced here for the solution of eq.(6) with the Cholesky decomposition of matrix $H$ (if the Cholesky decomposition fails for some values of sigma of the RBF, one can always replace $\Gamma$ with the identity matrix or use the least-squares solution to (6) based on the SVD for a very limited number of pixels). Furthermore, the positivity contraint in eq.(9) was not enforced as was done in [4], [25] since this increased the computation time while not providing any advantage for the final classification accuracy.

As opposed to the original reweighted $L_1$-norm scheme proposed in [32], Weighted $L_1$-norm SUnSAL is convex so it converges straight to its global minimum. At some values of sigma (for the RBF kernel) many of the eigenvalues of $A^T A$ may become almost zero resulting in a rather large flat region in the performance surface, but this performace surface does not become concave as in [32].

To finish the classification task, the reconstruction residuals (from the quadratic term in (3)) for each class are summed up for the closest (top) M neighbors in an NxN spatial window around each central pixel and the final decision is made in favor of the class which presents the minimum residual sum. Closeness to the central pixel in each window is measured based on the complement of the cosine of the angle between the central pixel and each neighbor or the Euclidean distance between the two vectors. As opposed to the pre-lowpass filtering procedure method used in [4], such a post-processing approach is more selective, but it does not force all the pixels in the selected neighborhood to be unmixed with the same support, as (K)SOMP, (K)SSP [2] or (K)CL-SUnSAL [26].

## IV. RESULTS

### A. Indian Pines image

The first results obtained on this image have bee initially presented in [41]; this section updates the results for all algorithms (now run on the same test sets) and adds the results obtained by using composite kernels (CK's).

The image used in our experiments is the Indian Pines scene from the Airborne Visible/Infrared Imaging Spectrometer (AVIRIS) [38] with a spatial resolution per pixel of 20 meters. It has 220 bands across the spectral range from 0.2 to 2.4μm, but in the experiments 20 water absorption bands [2], [4] have been removed (bands 104-108, 150-163, 220). There are 16 ground truth classes and 10% of the samples in each class were randomly chosen as training samples (see Appendix A, Table A1.) in each of our experiments as in [2], [4] for comparison purposes. Note that this method of sampling leads to undersampling of several classes (such as classes 1, 7, 9, 15, 16) and other works [3], [36] use only a reduced set of 9 classes for which an adequate number of training samples can be provided.

The classification performance of each of the 16 classes averaged over 20 random trials, the overall accuracy (OA), the average accuracy (AA) [2] and the kappa coefficient of agreement $\kappa$ [39], [40] are shown in Tables 1- 6 for various algorithms (kernelized (KOMP / KSP) or not kernelized (OMP / SP), with / without spatial processing (SOMP / SSP versus OMP / SP)).

All kernelized algorithms used an RBF kernel for this image and all hyperspectral pixel vectors (both dictionary and test pixels) were normalized to unit norm for the sparse unmixing methods. The neighborhood window size was set to N=9 and M= 45 to 65 nearest neighbors were used for each central pixel in the spatial post-processing phase for Tables 1, 2, 5 and 6. These choices were made based on cross-validation in order to maximize classification accuracy.

The procedure in **Algorithm 2** was used with $\Gamma = I$ and then with $\Gamma = \Gamma$, as shown in columns **A** and **B** of Table 1. Both Weighted $L_1$-norm SUnSAL and KSOMP / KSSP have been kernelized with $\sigma = 250$ and the cosine of the angle between vectors was used to select the top M closest neighbors of each central pixel. For both columns **A** and **B**, the weights in $\Gamma$ were designed based on the cosine of the angle between the current pixel $y$ and each column of the dictionary matrix $A$ in the kernel induced high-dimensional Hilbert space.

The standard deviation of the OA is 0.27% for both columns A and B in Table 1. It was observed that our proposed algorithm was most effective in improving classification results when the rank of the Hessian $A^T A$ was considerably lower than its dimension (e.g. in the original data space where the rank is 200 for a dimension of 1043x1043), as the proposed algorithm uses 'locality' to compensate for the rank deficiency of $H$ (see Table 5). We see from Table 1 that our method performs better than the SVM with spectral-spatial composite kernels, but about the same as KSOMP and KSSP. It should be noted that the OA result given by KSOMP is rather unstable: depending on the training set used to perform validation, the final OA result can be from 0.5 – 1.5% lower than what is shown in Table 1. More results given by other classifiers can be found in [2], [4] for comparison purposes (it should be noted that these two references only used the spatial mean in $x_i^s$ for the CK [15] paradigm, see eq. (10)) . It is worth noting that both KSOMP, KSSP and the original L$_1$-norm SUnSAL can sometimes achieve higher OA, AA and $\kappa$ results than shown if they are used with the parameters detected by validation by the proposed method, Weighted L$_1$-norm SUnSAL (at least for this image, e.g. Table 6) ; however, these algorithms are not always capable of detecting these correct parameters by validation on their own.

---

[2] Note that the quadratic term in (3) will always be convex as long as a positive definite kernel is used.

We see that even when run in the original data space (column 4 of Table 5), our algorithm is close in performance to SVM-CK [15] in terms of OA and $\kappa$. When kernelized (column B of Table 1), the performance advantage given by the introduction of matrix $\Gamma$ in eq.(3) is much smaller (negligible) than in the original data space, as the rank of the Hessian $A^T A$ is 1043, equal to the matrix dimension (often with a low condition number). The proposed method still obtains higher OA values by about 0.5-0.6% than the original SUnSAL alg. for dictionaries $A$ where the Hessian is rank deficient, but these cases do not alter the final average over 20 random trials.

Joint spatial processing is known to drastically reduce the classification accuracy for classes with small spatial extent (e.g. class 9). Unlike the other classification methods included in Tables 1-6 where only the final spatially smoothed result is available, our method allows the inspection of the classification results based on simple unmixing only (before spatial postprocessing is applied), where the detection accuracies for classes with a very small spatial extent are much higher (e.g. compare the results for Table 5, class 9, columns 2 and 4 where 67.78%>>0%.).

Table 2 shows how the adaptively weighted L$_1$-norm SUnSAL algorithm compares with algorithms approximating L$_0$-norm regularization to induce sparsity, namely SOMP and SSP [2], when only spatial smoothing is used in the original data space. It is observed that SOMP and SSP give higher results in terms of OA, AA and coefficient of agreement $\kappa$ than the proposed algorithm, but this advantage vanishes when both spatial smoothing and kernelization are used together, as shown in the last three columns of Table 1.

Tables 3 and 4 compare the performance of our algorithm (without spatial postprocessing) with the L$_0$-norm 'regularized' OMP and SP in the original data space and after kernelization, respectively. We can see that our algorithm performs better than (K)OMP and (K)SP in the OA, AA and $\kappa$ coefficient values, the differences in performance being smaller for the kernelized algorithms than in the original data space. Even in kernel space, however, the KOMP and KSP algorithms tend to be caught in some unfavorable local minimum, as evidenced by the results for classes 1, 4, 9, 10 and 15.

Table 5 compares the performance of the standard SUnSAL algorithm ($\Gamma = I$) with the adaptively weighted L$_1$-norm SUnSAL algorithm ($\Gamma = \Gamma$) in the original data space, with and without spatial postprocessing. Our method, which benefits from the introduction of $\Gamma$ in eq. (3), is superior in terms of OA, AA and kappa coefficient of agreement.

Finally, Table 6 shows that our algorithm outperforms both KSOMP-CK and KSSP-CK when the joint sparsity model (JSM) [2] and the composite kernel (CK) [15] paradigm (see eq. (10)) are combined, as both KSOMP-CK and

$$k\left(x_i, x_j\right) = \mu \cdot k_s\left(x_i^s, x_j^s\right) + (1-\mu) \cdot k_w\left(x_i^w, x_j^w\right) (10)$$

KSSP-CK seem to be caught in a local unfavorable minimum for the final OA results (a full 9x9 window was used for the extraction of the spatial mean and standard deviation used in the spatial vector $x_i^s$ for the CK in this table).

### B. Pavia University image

The image used in the experiments of this section is the urban image of the 'University of Pavia' acquired by the Reflective Optics System Imaging Spectrometer (ROSIS). The ROSIS sensor generated 115 spectral bands ranging from 0.43 to 0.86 μm with a spatial resolution of 1.3m per pixel (of which the 12 noisiest bands were removed [2]).

This section presents a reduced set of results in Tables 7 and 8 (for a 5x5 and 7x7 spatial postprocessing window, respectively) where the performance of the Weighted L$_1$-norm SUnSAL algorithm is compared to the performance of (KS)OMP and (KS)SP with the same (original) training and test sets as in [2] (see Appendix A, Table A2.) for a few cases of interest, namely in the original data space, in kernel space (spectral kernel only) and in the weighted sum of spatial and spectral components CK space, as in eq.(10). All hyperspectral pixel vectors (both dictionary and test pixels) were normalized to unit norm for the sparse unmixing methods and in all kernelized algorithms an RBF kernel was used, just as for the Indian Pines image in the previous section.

For this image an adaptive number of pixels was employed for the spatial postprocessing phase by cutting off the neighbors which were added to the final reconstruction residual sum when the sum of their distances to the central pixel of the neighborhood exceeded a certain slope. This adaptive neighborhood scheme rarely gave an advantage for the final classification OA in Table 7 versus using the entire neighborhood of pixels as in [2], but it allowed to obtain lower simulation run times as fewer pixels were unmixed overall. Even with this improvement, the OA results for KSOMP / KSSP initially published in [2] for this image could not be exactly reproduced, at least as long as the RBF kernel was used.

**Table 1**. Weighted L$_1$-norm SUnSAL algorithm compared with other top of the line classification methods, with spatial postprocessing (Indian Pines image).

| Class | SVM-CK | A kernelized Original L$_1$-norm SUnSAL | B kernelized Weighted L$_1$-norm SUnSAL | K-SOMP | K-SSP |
|---|---|---|---|---|---|
| 1 | 95.63 | 94.90 | 94.90 | 94.90 | 94.90 |
| 2 | 95.76 | 95.92 | 95.92 | 96.04 | 96.02 |
| 3 | 97.18 | 96.39 | 96.39 | 96.51 | 96.49 |
| 4 | 93.00 | 95.29 | 95.29 | 95.38 | 95.40 |
| 5 | 96.24 | 97.06 | 97.06 | 97.06 | 97.06 |
| 6 | 98.76 | 98.72 | 98.72 | 98.72 | 98.72 |
| 7 | 84.35 | 89.35 | 89.35 | 89.35 | 89.35 |
| 8 | 99.31 | 100 | 100 | 100 | 100 |
| 9 | 87.50 | 38.06 | 38.06 | 38.06 | 38.06 |
| 10 | 94.12 | 96.21 | 96.21 | 96.25 | 96.25 |
| 11 | 96.22 | 98.97 | 98.97 | 98.94 | 98.94 |
| 12 | 95.34 | 96.28 | 96.28 | 96.29 | 96.30 |
| 13 | 99.58 | 98.53 | 98.53 | 98.53 | 98.53 |
| 14 | 98.44 | 100 | 100 | 100 | 100 |
| 15 | 94.90 | 97.98 | 97.98 | 97.98 | 97.98 |
| 16 | 96.65 | 96.82 | 96.82 | 96.82 | 96.82 |
| **OA** | 96.49 | 97.68 | 97.68 | 97.71 | 97.70 |
| **AA** | 95.18 | 93.16 | 93.16 | 93.18 | 93.18 |
| **κ** | 0.960 | 0.974 | 0.974 | 0.974 | 0.974 |

**Table 2**. Weighted L$_1$-norm SUnSAL algorithm compared with other classification methods in the original data space, with spatial processing (Indian Pines image).

| Class | Original data space, with spatial postprocessing | | |
|---|---|---|---|
| | SOMP | SSP | Weighted L$_1$-norm SUnSAL |
| 1 | 89.38 | 91.46 | 84.38 |
| 2 | 95.53 | 94.76 | 91.58 |
| 3 | 92.42 | 93.21 | 87.78 |
| 4 | 92.31 | 92.90 | 94.02 |
| 5 | 93.59 | 95.82 | 94.99 |
| 6 | 99.74 | 99.25 | 99.91 |
| 7 | 70.65 | 72.61 | 72.39 |
| 8 | 100 | 99.99 | 100 |
| 9 | 3.33 | 31.67 | 0 |
| 10 | 87.76 | 95.12 | 91.71 |
| 11 | 98.01 | 98.25 | 99.23 |
| 12 | 92.64 | 96.94 | 97.84 |
| 13 | 99.42 | 99.21 | 99.18 |
| 14 | 99.61 | 99.85 | 99.88 |
| 15 | 96.40 | 92.02 | 79.66 |
| 16 | 98.18 | 96.06 | 98.00 |
| **OA** | 95.70 | 96.65 | 95.25 |
| **AA** | 88.06 | 90.57 | 86.91 |
| **κ** | 0.951 | 0.962 | 0.945 |

**Table 3**. Weighted L$_1$-norm SUnSAL algorithm compared with other classification methods in the original data space, no spatial processing (Indian Pines image).

| Class | Original data space, without spatial postprocessing | | |
|---|---|---|---|
| | OMP | SP | Weighted L$_1$-norm SUnSAL |
| 1 | 23.65 | 64.69 | 62.40 |
| 2 | 79.55 | 66.71 | 78.28 |
| 3 | 43.37 | 63.17 | 67.71 |
| 4 | 23.76 | 50.55 | 61.86 |
| 5 | 77.32 | 90.21 | 92.04 |
| 6 | 97.96 | 96.54 | 97.92 |
| 7 | 0 | 72.83 | 63.91 |
| 8 | 99.80 | 98.57 | 99.48 |
| 9 | 2.22 | 51.94 | 67.78 |
| 10 | 26.49 | 73.00 | 74.48 |
| 11 | 83.82 | 81.93 | 87.69 |
| 12 | 45.53 | 60.98 | 77.64 |
| 13 | 99.21 | 98.79 | 99.50 |
| 14 | 98.53 | 96.14 | 97.51 |
| 15 | 40.44 | 45.03 | 57.08 |
| 16 | 93.82 | 93.24 | 91.65 |
| **OA** | 72.44 | 78.47 | 84.02 |
| **AA** | 58.47 | 75.27 | 79.81 |
| **κ** | 0.679 | 0.753 | 0.817 |

**Table 4**. Kernelized Weighted L$_1$-norm SUnSAL algorithm compared with other classification algorithms, no spatial processing (Indian Pines image).

| Class | Kernelized, without spatial postprocessing | | |
|---|---|---|---|
| | KOMP | KSP | Weighted L$_1$-norm SUnSAL |
| 1 | 62.60 | 64.48 | 71.67 |
| 2 | 78.77 | 80.89 | 84.96 |
| 3 | 70.09 | 72.73 | 75.27 |
| 4 | 67.93 | 66.71 | 71.24 |
| 5 | 91.90 | 92.36 | 92.96 |
| 6 | 96.70 | 96.70 | 96.79 |
| 7 | 70.22 | 75.22 | 77.61 |
| 8 | 98.89 | 98.90 | 99.17 |
| 9 | 67.22 | 70.00 | 76.94 |
| 10 | 69.99 | 73.23 | 81.61 |
| 11 | 86.34 | 86.09 | 88.23 |
| 12 | 79.80 | 81.22 | 84.55 |
| 13 | 98.89 | 99.03 | 99.37 |
| 14 | 96.53 | 96.49 | 96.28 |
| 15 | 55.39 | 56.58 | 62.75 |
| 16 | 91.71 | 91.53 | 91.06 |
| **OA** | 83.50 | 84.40 | 87.07 |
| **AA** | 80.19 | 81.39 | 84.40 |
| **κ** | 0.811 | 0.822 | 0.852 |

**Table 5.** Original L$_1$-norm SUnSAL compared to the Weighted L$_1$-norm SUnSAL algorithm in the original data space (Indian Pines image).

| Class | Before spatial post-processing | | After spatial postprocessing | |
|---|---|---|---|---|
| | Original SUnSAL | Weighted L$_1$-norm SUnSAL | Original SUnSAL | Weighted L$_1$-norm SUnSAL |
| 1 | 45.31 | 62.40 | 42.08 | 84.38 |
| 2 | 76.16 | 78.28 | 88.23 | 91.58 |
| 3 | 44.80 | 67.71 | 67.63 | 87.78 |
| 4 | 36.43 | 61.86 | 70.38 | 94.02 |
| 5 | 84.83 | 92.04 | 89.74 | 94.99 |
| 6 | 97.65 | 97.92 | 99.99 | 99.91 |
| 7 | 6.30 | 63.91 | 0.00 | 72.39 |
| 8 | 99.63 | 99.48 | 100 | 100 |
| 9 | 14.72 | 67.78 | 0.00 | 0 |
| 10 | 35.71 | 74.48 | 54.00 | 91.71 |
| 11 | 87.78 | 87.69 | 99.68 | 99.23 |
| 12 | 54.38 | 77.64 | 90.16 | 97.84 |
| 13 | 99.32 | 99.50 | 99.37 | 99.18 |
| 14 | 97.64 | 97.51 | 99.92 | 99.88 |
| 15 | 48.33 | 57.08 | 71.46 | 79.66 |
| 16 | 91.82 | 91.65 | 98.29 | 98.00 |
| OA | 75.34 | 84.02 | 87.83 | 95.25 |
| AA | 63.80 | 79.81 | 73.18 | 86.91 |
| $\kappa$ | 0.713 | 0.817 | 0.859 | 0.945 |

**Table 6.** Composite kernel Weighted L$_1$-norm SUnSAL algorithm compared with other classification methods (Indian Pines image).

| Class | Kernelized, with spatial postprocessing | | | |
|---|---|---|---|---|
| | KSOMP-CK | KSSP-CK | Original SUnSAL-CK | Weighted L$_1$-norm SUnSAL-CK |
| 1 | 96.88 | 96.56 | 95.73 | 96.04 |
| 2 | 96.51 | 96.01 | 98.00 | 97.83 |
| 3 | 98.25 | 98.03 | 98.81 | 99.11 |
| 4 | 97.38 | 97.14 | 96.02 | 95.81 |
| 5 | 96.29 | 96.30 | 96.73 | 96.95 |
| 6 | 99.75 | 99.81 | 99.75 | 99.69 |
| 7 | 87.17 | 88.70 | 87.17 | 93.26 |
| 8 | 100 | 100 | 100 | 100 |
| 9 | 17.78 | 16.11 | 19.44 | 33.89 |
| 10 | 95.88 | 95.86 | 95.74 | 95.72 |
| 11 | 99.01 | 98.96 | 99.14 | 99.06 |
| 12 | 97.66 | 97.56 | 97.90 | 97.46 |
| 13 | 99.13 | 99.05 | 98.18 | 99.16 |
| 14 | 99.66 | 99.72 | 99.75 | 99.76 |
| 15 | 97.89 | 97.63 | 99.09 | 98.86 |
| 16 | 96.76 | 96.65 | 95.71 | 97.76 |
| OA | 97.99 | 97.88 | 98.31 | 98.32 |
| AA | 92.25 | 92.13 | 92.32 | 93.77 |
| $\kappa$ | 0.977 | 0.976 | 0.981 | 0.981 |

**Table 7.** Weighted L$_1$-norm SUnSAL (Proposed Algorithm) compared with other classification methods for the Pavia University image (5x5 spatial postprocessing window, 5x5 window for $x_i^s$ in col.'s 9-13 for the CK).

| | 1 | 2 | 3 | 4 | 5 | 6 | 7 | 8 | 9 | 10 | 11 | 12 | 13 |
|---|---|---|---|---|---|---|---|---|---|---|---|---|---|
| Class | Original data space, with joint spatial processing | | | | Kernelized, with joint spatial processing | | | | | Kernelized (CK) with joint spatial processing | | | |
| | SOMP | SSP | Original SUnSAL | Proposed algorithm | KSOMP | KSSP | Original SUnSAL | Proposed algorithm | SVM-CK | K-SOMP-CK | K-SSP-CK | Original SUnSAL | Proposed algorithm |
| 1 | 86.02 | 87.69 | 88.15 | 92.48 | 96.27 | 91.01 | 94.08 | 93.10 | 86.90 | 92.42 | 93.12 | 96.37 | 94.40 |
| 2 | 63.47 | 69.05 | 65.40 | 74.41 | 66.36 | 69.17 | 67.87 | 76.02 | 83.42 | 96.57 | 97.06 | 90.52 | 97.38 |
| 3 | 93.39 | 87.49 | 92.45 | 79.56 | 84.52 | 84.96 | 75.98 | 80.33 | 71.24 | 85.01 | 85.62 | 85.73 | 82.76 |
| 4 | 97.87 | 95.09 | 97.91 | 99.42 | 98.90 | 96.43 | 98.87 | 99.04 | 97.46 | 99.00 | 99.31 | 99.28 | 99.73 |
| 5 | 100 | 100 | 100 | 100 | 100 | 100 | 100 | 100 | 99.37 | 99.73 | 100 | 99.73 | 100 |
| 6 | 99.23 | 99.06 | 99.89 | 94.18 | 98.82 | 99.02 | 98.71 | 91.21 | 95.47 | 81.87 | 84.62 | 93.66 | 84.19 |
| 7 | 97.15 | 97.15 | 97.86 | 96.84 | 95.62 | 96.84 | 94.29 | 94.80 | 91.74 | 97.86 | 98.37 | 96.64 | 95.41 |
| 8 | 61.12 | 93.91 | 61.53 | 97.95 | 92.87 | 91.71 | 98.16 | 97.86 | 95.81 | 95.27 | 94.59 | 94.95 | 95.75 |
| 9 | 98.49 | 97.61 | 99.87 | 99.75 | 99.50 | 98.49 | 99.62 | 99.87 | 96.35 | 99.25 | 98.99 | 99.25 | 99.25 |
| OA | 77.32 | 82.36 | 78.64 | 85.32 | 82.52 | 82.74 | 82.87 | 85.76 | 87.76 | 93.95 | 94.61 | 93.17 | 94.83 |
| AA | 88.53 | 91.89 | 89.23 | 92.73 | 92.54 | 91.96 | 91.95 | 92.47 | 90.86 | 94.11 | 94.63 | 95.13 | 94.32 |
| $\kappa$ | 0.716 | 0.776 | 0.732 | 0.812 | 0.779 | 0.781 | 0.783 | 0.817 | 0.841 | 0.918 | 0.927 | 0.909 | 0.930 |

**Table 8.** Weighted L$_1$-norm SUnSAL (Proposed Algorithm) compared with other classification methods for the Pavia University image (7x7 spatial postprocessing window, 3x3 window for $x_i^s$ in col.'s 9-13 for the CK).

| | 1 | 2 | 3 | 4 | 5 | 6 | 7 | 8 | 9 | 10 | 11 | 12 | 13 |
|---|---|---|---|---|---|---|---|---|---|---|---|---|---|
| Class | Original data space, with joint spatial processing | | | | Kernelized, with joint spatial processing | | | | | Kernelized (CK) with joint spatial processing | | | |
| | SOMP | SSP | Original SUnSAL | Proposed algorithm | KSOMP | KSSP | Original SUnSAL | Proposed algorithm | SVM-CK | KSOMP-CK | KSSP-CK | Original SUnSAL | Proposed algorithm |
| 1 | 86.28 | 86.47 | 89.40 | 92.78 | 97.32 | 87.01 | 94.19 | 94.21 | 90.97 | 93.43 | 96.84 | 95.65 | 95.75 |
| 2 | 64.29 | 70.11 | 65.66 | 76.75 | 67.57 | 70.79 | 68.56 | 78.29 | 81.37 | 96.48 | 95.10 | 91.70 | 96.93 |
| 3 | 93.28 | 86.12 | 92.34 | 79.95 | 84.68 | 84.08 | 76.80 | 80.33 | 77.19 | 88.04 | 87.38 | 89.20 | 86.78 |
| 4 | 96.84 | 94.47 | 98.70 | 98.73 | 98.15 | 96.67 | 97.97 | 97.87 | 98.11 | 99.38 | 99.35 | 99.52 | 99.69 |
| 5 | 100 | 100 | 99.91 | 100 | 100 | 100 | 99.73 | 100 | 99.91 | 100 | 100 | 100 | 100 |
| 6 | 99.72 | 99.54 | 99.98 | 95.84 | 99.06 | 98.99 | 98.78 | 93.83 | 97.53 | 88.34 | 94.60 | 97.44 | 96.48 |
| 7 | 96.84 | 96.33 | 98.27 | 97.86 | 96.53 | 96.53 | 97.25 | 96.74 | 90.93 | 98.06 | 98.67 | 98.98 | 97.66 |
| 8 | 64.42 | 95.57 | 63.32 | 97.74 | 94.59 | 92.98 | 98.66 | 97.53 | 95.66 | 95.21 | 94.80 | 94.41 | 97.27 |
| 9 | 92.83 | 92.70 | 99.37 | 99.25 | 97.74 | 95.22 | 98.99 | 98.99 | 97.86 | 94.59 | 93.84 | 95.85 | 93.59 |
| OA | 77.86 | 82.62 | 79.16 | 86.58 | 83.35 | 82.85 | 83.28 | 87.18 | 88.04 | 94.89 | 95.45 | 94.15 | 96.50 |
| AA | 88.28 | 91.26 | 89.66 | 93.21 | 92.85 | 91.36 | 92.33 | 93.09 | 92.17 | 94.62 | 95.62 | 95.86 | 96.02 |
| κ | 0.722 | 0.779 | 0.738 | 0.827 | 0.789 | 0.782 | 0.788 | 0.834 | 0.845 | 0.931 | 0.939 | 0.922 | 0.953 |

We observe that the Weighted L$_1$-norm SUnSAL algorithm gives higher OA and kappa results in each category. Contrary to the Indian Pines image case, the classification accuracy maxima obtained for the Pavia University image occur at values of the parameter σ for the RBF kernel at which the rank of the Hessian $A^T A$ is considerably lower than its dimension of 3921, a fact which allows the proposed algorithm to work efficiently.

On the other hand, KSOMP and KSSP need these maxima to occur where the Hessian $A^T A$ is full rank with a low condition number in order for them to work in a stable manner, a condition which was present for the Indian Pines image but not for the Pavia University image.

The results obtained for the original SUnSAL algorithm ($\Gamma = I$) are also shown in Table 7 in order to be able to quantify the advantage of the introduction of matrix $\Gamma$ in eq.(3). It is worth noting that the performance of KSOMP is also rather sensitive to the training set used as columns of the dictionary $A$. A simple change of the dictionary columns for class 2, Meadows, for example, can lead to drops in performance of the order of 4-5% if the parameters obtained by validation with the original training dictionary are used.

Just as for the Indian Pines image in Table 4, kernelized Weighted L$_1$-norm SUnSAL without spatial postprocessing gives higher OA results than KSP and KOMP (77.31% compared to 76.62% and 72.53%, respectively), showing that the quality of the unmixing solution is higher for the proposed method.

Also, the authors of [2] apparently did not realize that the joint sparsity model (JSM) they introduced in their algorithms is to a certain extent complementary to the spatial-spectral composite kernel paradigm and did not use the two together. By combining the JSM and the spatial-spectral CK models, OA results higher by about 11-12% than those shown in [2] for KOMP-CK and KSP-CK were obtained for KSOMP-CK and KSSP-CK as shown in columns 10 and 11 of Table 7.

The lack of randomness in the selection of signatures for the subdictionary for class 2, Meadows, is apparently responsible for the low classification results obtained on the Pavia University image [2]. To confirm this hypothesis, I have changed the subdictionary for class 2, Meadows, with a randomly selected dictionary from the entire class. The significantly improved classification results are shown for several algorithms in Fig. 2.

It is also interesting to point out that the proposed algorithm gives relatively stable results when run over the parameter σ of the RBF kernel function, contrary to the behaviour of KSOMP, KSSP or that of the original SUnSAL algorithm, as shown in Fig. 2. This characteristic may be very important when there is an insufficient number of training / test samples for validation and errors in the estimation of σ may occur. KSOMP and KSSP perform rather well at the left of the the graph where the rank of the Hessian matrix $A^T A$ is full, then their performance drops relatively fast as the rank of the Hessian decreases; a similar conclusion can be drawn from the behaviour of the original SUnSAL algorithm. On the other hand, the proposed method, Weighted $L_1$-norm SUnSAL, gives a relatively high classification OA value across all sigma values, reaching its peak at about σ=17e+3, just slightly higher than the peaks of KSOMP / KSSP on the left of the graph.

Last, but not least, the adaptive neighborhood selection scheme used for this image allowed the use of a variable subset of a 7x7 neighborhood postprocessing window instead of just a full 5x5 window as in [2]. This enabled all the algorithms to base their spatial processing only on those pixels in the 7x7 neighborhood most closely related to the central pixel of interest. The improved results are shown in Table 8, where the proposed algorithm again gives higher classification results than the competition (a full 3x3 window was used for the extraction of the spatial mean and standard deviation (per band) used in the spatial vector $\boldsymbol{x}_i^s$ for the CK in col.'s 9-13 of this table).

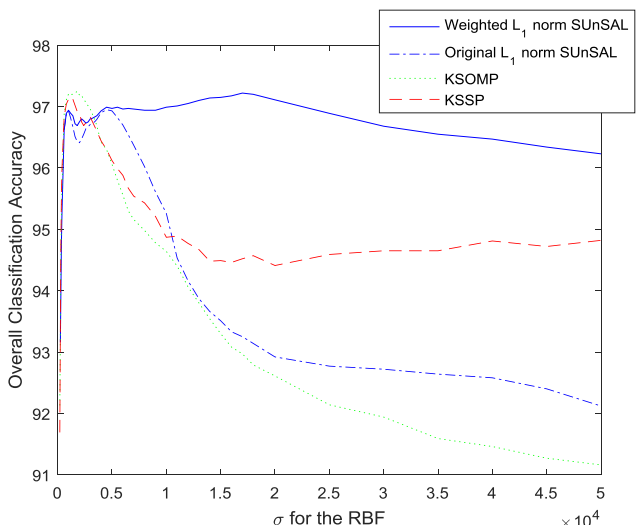

**Fig. 2** Performance of the proposed method, Weighted $L_1$-norm SUnSAL, compared to that of KSOMP, KSSP and original SUnSAL versus σ for the RBF kernel function (University of Pavia image, random sub-dictionary used for class 2 Meadows, 5x5 spatial postprocessing window).

## V.  CONCLUSION

In this paper I have introduced a new sparse regression method for hyperspectral image classification based upon the SUnSAL algorithm [24], [25], which adaptively weighs the level of sparsity applied to each pixel. As opposed to the original reweighted $L_1$-norm scheme in [32], Weighted $L_1$-norm SUnSAL is convex so it converges straight to its global minimum without running the risk of getting stuck in a local, suboptimal minimum. Unlike other heuristic adaptive $L_1$-norm schemes, the proposed method works well because its adaptive weights are based upon the actual physical separation of hyperspectral signatures.

Whereas the SVM is a discriminative classification method, the proposed algorithm can be seen as a generative model similar to (KS)OMP and (KS)SP, where the subspaces representing different classes compete with each other during the pixel unmixing process, leading to a vector of unmixing coefficients which has only a few non-zero, representative coefficients.

From the results obtained on the Indian Pines image, we see that the proposed algorithm does not get caught in local unfavorable minima, as is the case for (KS)OMP and (KS)SP which are not convex optimization problems. Because the rank of the Hessian $A^T A$ of the data reconstruction fidelity term is full in kernel space where the classification maxima for the Indian Pines image occur, the final classification results for the proposed method and those of KSOMP and KSSP are rather close in the last three columns of Table 1. However, the classification maxima for the University of Pavia image in Tables 7 and 8 occur where the rank of the Hessian $A^T A$ is lower than its dimension, a fact which allows the proposed algorithm to work efficiently and obtain higher results than all the other algorithms. It is expected that by using more sophisticated CK's such as those based on EMAP's, as in [42], [43], classification accuracies in the vicinity of 98% could be obtained for the University of Pavia image with the original training dictionary and the proposed method. The experiments conducted so far show that the adaptively Weighted $L_1$-norm SUnSAL algorithm will not adversely affect the final classification accuracy when the classification OA maxima occur where the rank of the Hessian $A^T A$ is full in kernel space, as in Table 1, but will provide a clear advantage when these classification maxima occur where the rank of the Hessian $A^T A$ is much smaller than its dimension, as in Tables 7 and 8.

This paper has shown it is possible to design a convex, adaptive, $L_1$-norm regularized sparse regression method whose results are superior to those of state of the art algorithms such as SVM-CK and reach or surpass the performance of non-convex methods such as KSOMP and KSSP, while at the same time providing full pixel unmixing results, an additional feature not available from the other algorithms.

As images with large amounts of fine detail (such as the Pavia University image) become increasingly more common as hyperspectral sensors evolve, it is expected that the algorithm presented in this paper will also gain in relative importance versus existing non-convex counterparts. Other applications include convex hyperspectral image unmixing, target detection and possibly more areas in the domain of compressive sensing.

## ACKNOWLEDGEMENT

The author would like to thank Professor Paolo Gamba, the University of Pavia and the ROSIS project for providing the data for the University of Pavia hyperspectral image used in the experiments in Tables 7 and 8. The libSVM software from http://www.csie.ntu.edu.tw/~cjlin/ was used to obtain the results for SVM-CK [15].

## VI. APPENDIX A

This appendix presents the number of training and test samples used for the Indian Pines and Pavia University images in Section IV.

TABLE A1. Training / Test Samples for the Indian Pines image [38].

| No. | Class Name | Train | Test |
|---|---|---|---|
| 1 | Alfalfa | 6 | 48 |
| 2 | Corn-notill | 144 | 1290 |
| 3 | Corn-min | 84 | 750 |
| 4 | Corn | 24 | 210 |
| 5 | Grass/Pasture | 50 | 447 |
| 6 | Grass/Trees | 75 | 672 |
| 7 | Grass/Pasture mowed | 3 | 23 |
| 8 | Hay-windrowed | 49 | 440 |
| 9 | Oats | 2 | 18 |
| 10 | Soybeans-notill | 97 | 871 |
| 11 | Soybeans-min | 247 | 2221 |
| 12 | Soybean-clean | 62 | 552 |
| 13 | Wheat | 22 | 190 |
| 14 | Woods | 130 | 1164 |
| 15 | Building-Grass-Trees-Drives | 38 | 342 |
| 16 | Stone-steel Towers | 10 | 85 |
|  | Total | 1043 | 9323 |

TABLE A2. Training / Test Samples for the Pavia University image.

| No. | Class Name | Train | Test |
|---|---|---|---|
| 1 | Asphalt | 548 | 6304 |
| 2 | Meadows | 540 | 18146 |
| 3 | Gravel | 392 | 1815 |
| 4 | Trees | 524 | 2912 |
| 5 | Metal Sheets | 265 | 1113 |
| 6 | Bare Soil | 532 | 4572 |
| 7 | Bitumen | 375 | 981 |
| 8 | Bricks | 514 | 3364 |
| 9 | Shadows | 231 | 795 |
|  | Total | 3921 | 40002 |